\documentclass[a4paper,11pt]{article}

\usepackage{amsmath,amsthm,amssymb,amsfonts}
\usepackage{mathtools}
\usepackage{graphicx}
\usepackage{float}
\usepackage{tikz}
\usepackage{tikz-cd}
\usepackage{xcolor}
\usepackage{textcomp}
\usepackage{booktabs}
\usepackage{multirow}
\usepackage[title]{appendix}
\usepackage{hyperref}
\usepackage{listings}
\usepackage{algorithm}
\usepackage{algorithmicx}
\usepackage{algpseudocode}
\usepackage{mathrsfs}

\mathtoolsset{showonlyrefs,showmanualtags}

\newtheorem{theorem}{Theorem}[section]
\newtheorem{lemma}[theorem]{Lemma}

\theoremstyle{definition}
\newtheorem{definition}[theorem]{Definition}
\newtheorem{example}[theorem]{Example}
\theoremstyle{remark}

\DeclareMathOperator{\rkhs}{\textsf{RKHS}}
\DeclareMathOperator{\tr}{tr}
\DeclareMathOperator{\Li}{Li}

\title{Commutators of the Unilateral Shift and Adjoint for Reproducing Kernel Hilbert Spaces on the Disk}
\author{Nathan Parker\thanks{The author was supported by an EPSRC scholarship. Special thanks to Gordon Blower for his guidance and support.}}
\date{} % Leave blank or use \date{\today}

\newcommand{\keywords}[1]{\noindent\textbf{Keywords: }#1}
\newcommand{\classifications}[1]{\noindent\textbf{MSC (2020): }#1}

\begin{document}

\maketitle


\classifications{30H05, 47A07, 47B37}

\begin{abstract}
We generalise the result of Berger and Shaw \cite{berger2006intertwining} the trace formula for Hardy Hilbert space to a larger class of rotation invariant Hilbert function spaces on the unit disk. We also demonstrate many meaningful examples of these Hilbert spaces by computing the inner products. We also extend to a wider class than the unilateral shift, that is, weighted shifts under certain restrictions.
\end{abstract}

%%%%%%%%%%%%%%%%%%%%%%%  AUTHOR's MAIN TEXT BEGINS HERE:
%%%%%%%%%%%%%%%%%%%%%%%  Use amsart conventions. 
%%%%%%%%%%%%%%%%%%%%%%

\section{Introduction}
This paper proves an extension of the Berger-Shaw theorem regarding the trace formula for the shift and its adjoint.  Berger and Shaw dealt with the Hardy Hilbert space on the disk while we extend to a class of rotation invariant Hilbert function spaces on the disk, remarkably, all these trace formulas involve Dirichlet space. A recent summary of the historical progress made relating to the trace formula can be found in \cite{howe2012traces}. 

\section{Reproducing Kernel Hilbert Spaces on the Disk}

\begin{definition}[Reproducing kernel Hilbert space] A reproducing kernel Hilbert space $\rkhs$ on a domain $D\subseteq \mathbb{C}$ is a complex Hilbert space $H$ of functions on $D$ such that the maps of point evaluations $f \rightarrow f(z)$ are continuous linear functionals.  For all $z \in D$ there exists a unique $K_z \in H$ such that $\langle f , K_z \rangle = f(z)$.  Let $K(z , w) = \langle K_w , K_z \rangle_H$.  
\end{definition}

\begin{lemma} Let $H$ be a $\rkhs$: \begin{enumerate} 
\item{Suppose $z \rightarrow K_z$ is a weakly continuous map $D \rightarrow H$ then the function $(z , w) \rightarrow K(z, w)$ is continuous for all $z, w $.  }
\item{Suppose further \begin{equation} \int_\Gamma K(z , w) dz =0 \end{equation} for all contours $\Gamma$ in $D$ then $f(z)= \langle f , K_z \rangle$ is holomorphic on $D$.  }
\end{enumerate}
\begin{proof} \begin{enumerate}\item{By the weak continuity of $z \rightarrow K_z$, the map $z \rightarrow \langle f , K_z\rangle $ is continuous for all $f \in H$.  We let $f= K_w$ and we have $z \rightarrow  \langle K_w , K_z \rangle = K(z , w) $ continuous.  The same argument holds for $w$ and by symmetry we have joint continuity of the map. 
 }
\item{This is due to Morera's theorem. }
\end{enumerate}
\end{proof}
\end{lemma}

\begin{definition}[Rotation Invariance] For the open unit disc $\mathbb{D}$, say that  a reproducing kernel Hilbert space is rotation invariant if $R_\theta :f(z)\mapsto f(e^{i\theta}z)$ gives a linear isometric isomorphism on $H$.\end{definition}

An invertible isometry is a unitary, so $R_{- \theta} = R_\theta^\dagger$.

\begin{definition}[Dirichlet Space] The Dirichlet space $\mathcal{D}$ on the unit disk, $\mathbb{D}$ is the space of holomorphic functions such that for all $f \in \mathcal{D}$, we have \begin{equation}\displaystyle \int_{\mathbb{D}} |f'(z)|^2 dA(z)<\infty . \end{equation} The inner product is given by \begin{equation}\langle f , g \rangle_{\mathcal{D}} = f(0) \overline{g(0)} + \frac{1}{\pi} \displaystyle \int_{\mathbb{D}} f'(z) \overline{g'(z)} dA(z) . \end{equation}

The Dirichlet space $\mathcal{D}$ gives a $\rkhs$ on $\mathbb{D}$.  Let $\mathcal{D}_0$ be the closed linear subspace $\mathcal{D}_0=\{ f\in \mathcal{D}: f(0)=0\}$ of $\mathcal{D}$.  The orthonormal basis of $\mathcal{D}_0$ is given by $\left\{ \frac{1}{\sqrt{n} } z^n \right\}_{n=1}^\infty$. Let \begin{equation}f(z)=\sum_{j=0}^\infty a_nz^n,\qquad   g(z) =\sum_{n=0}^\infty b_nz^n ,\end{equation} then 
\begin{equation} \langle f,g\rangle =\sum_{n=1}^\infty (n+1) a_n\overline{b_n}.\end{equation}
\end{definition}

\noindent From here onwards we assume $\left( \alpha_n \right)_{n=0}^\infty$ is a sequence of positive real numbers.  

\begin{definition} For each sequence $\alpha=(\alpha_n)_{n=0}^\infty$ such that \begin{equation}\displaystyle \limsup_{n \rightarrow \infty} (\alpha_n)^{\frac{1}{n}} = 1 ,\end{equation}
 let $H_{\alpha}$ be the Hilbert space with elements power series $f,g\in H_{\alpha}$ given by \begin{equation}f(z)=\displaystyle \sum_{n=0}^\infty a_n z^n , \qquad g(z)=\displaystyle \sum_{n=0}^\infty b_n z^n\end{equation} with inner product given by \begin{equation}\langle f , g \rangle_{H_{\alpha}} =\displaystyle \sum_{n=0}^\infty a_n \overline{b_n} \alpha_n .\end{equation}   
\end{definition}

\section{Main Theorem}

\begin{theorem}

Let $\alpha=(\alpha_n)_{n=0}^\infty $ obey  \begin{equation}\displaystyle \lim_{n \rightarrow \infty} \frac{\alpha_{n+1}}{\alpha_{n}}=1\end{equation} and suppose $\alpha$ is concave or convex. 
 For \begin{equation}f(z)=\displaystyle \sum_{n=0}^\infty a_n z^n \end{equation} the following properties hold:

\begin{enumerate}
    \item $H_{\alpha}$ is a rotation invariant Hilbert space
    \item $H_{\alpha}$ has reproducing kernel \begin{equation}K_w (z) = \displaystyle \sum_{n=0}^\infty \frac{\overline{w}^n z^n}{\alpha_n}. \end{equation}
    \item Let $S$ be the unilateral shift $Sf(z)=zf(z)$. Then $S$ is a bounded linear operator on $H_{\alpha}$
    \item The adjoint shift $S^\dagger$ of \begin{equation}f(z)=\displaystyle \sum_{n=0}^\infty a_n z^n\end{equation} on $H_{\alpha}$ is given by \begin{equation}S^\dagger f(z)=\displaystyle \sum_{n=0}^\infty \frac{\alpha_{n+1}}{\alpha_n} a_{n+1} z^n .\end{equation}  
    \item The commutator of $S$ and $S^\dagger$ is trace class and for all polynomials $f$ and $g$  \begin{equation}\langle \langle f , g \rangle \rangle_{H_{a_n}} = \tr \left( g(S)^\dagger f(S) - f(S) g(S)^\dagger \right)=\frac{1}{\pi} \displaystyle \int_{\mathbb{D}} f'(z) \overline{g'(z)} dA(z)\end{equation}  
\end{enumerate}
    \begin{proof} 
    \begin{enumerate}
        \item We have \begin{equation} ||f(e^{i \theta} z )||_{H_\alpha}^2=\displaystyle\sum_{n=0}^\infty |a_n^2 e^{i n \theta}| \alpha_n= ||f||_{H_\alpha}^2.\end{equation}  
        \item We have \begin{equation}f(w)=\displaystyle \sum_{n=0}^\infty a_n w^n= \displaystyle \sum_{n=0}^\infty  \alpha_n a_n \frac{w^n}{\alpha_n}= \langle f , K_w  \rangle_{H_{\alpha}}.\end{equation}  Hence \begin{equation}K_w(z)= \displaystyle \sum_{n=0}^\infty \frac{\overline{w}^n z^n}{\alpha_n}.\end{equation}  
        \item Firstly, we have $S(\lambda f +g ) = z ( \lambda f + g )= \lambda zf + zg = \lambda Sf + Sg$ hence $S$ is linear.  Now we have \begin{equation}||f||_{H_{\alpha}}^2=\displaystyle \sum_{n=0}^\infty |a_n^2| \alpha_n\end{equation} and \begin{equation}||Sf||_{H_{\alpha}}^2=\displaystyle \sum_{n=0}^\infty |a_n^2| \alpha_{n+1}= \displaystyle \sum_{n=0}^\infty |a_n^2| \alpha_n \frac{\alpha_{n+1}}{\alpha_n}\end{equation} hence since we know the series obtained by \begin{equation}\displaystyle \sup_{n \in \mathbb{N}_0}\left\{\frac{\alpha_{n+1}}{\alpha_n}\right\} \displaystyle \sum_{n=0}^\infty |a_n^2| \alpha_n \end{equation} is convergent and we have \begin{equation} \displaystyle \sup_{n \in \mathbb{N}_0}\left\{\frac{\alpha_{n+1}}{\alpha_n}\right\} \displaystyle \sum_{n=0}^\infty |a_n^2| \alpha_n   \geq  \displaystyle \sum_{n=0}^\infty |a_n^2| \alpha_n \frac{\alpha_{n+1}}{\alpha_n} \end{equation} We must have \begin{equation}||Sf||_{H_{\alpha}}^2 \leq \displaystyle \sup_{n \in \mathbb{N}_0}\left\{\frac{\alpha_{n+1}}{\alpha_n}\right\} ||f||_{H_{\alpha}}^2.\end{equation}  Hence $S$ is a bounded operator.  
        \item The adjoint shift must satisfy $\langle Sf, g \rangle_{H_{\alpha}} = \langle f , S^\dagger g \rangle_{H_{\alpha}}$.  We have \begin{equation}\langle Sf , g \rangle_{H_{\alpha}}= \displaystyle \sum_{n=1}^\infty a_{n-1} \overline{b_n} \alpha_n ,   \langle f , S^\dagger g \rangle_{H_{\alpha}}= \displaystyle \sum_{n=0}^\infty a_n \overline{b_{n+1}} \alpha_{n+1}\end{equation} and these are equal by change of indices hence the operator described must be the adjoint shift.  
        \item We first consider $f(z)=g(z)=z^m$ and consider the operation on elements of the orthogonal basis $\left\{ z^n \right\}_{n=0}^\infty $.  We have: \begin{equation}\begin{tikzcd} z^n \arrow[rightarrow]{r}{S^m} &  z^{n+m}\arrow[rightarrow]{r}{(S^\dagger)^m} & \frac{\alpha_{n+m}}{\alpha_n} z^n\end{tikzcd}.  \end{equation}
    Also: \begin{equation}\begin{tikzcd} z^n \arrow[rightarrow]{r}{(S^\dagger)^m} &  \frac{n-m+ \frac{1}{2} + |n-m+\frac{1}{2}|}{2(n-m+\frac{1}{2})} \frac{\alpha_n}{\alpha_{n-m}}z^{n-m}\arrow[rightarrow]{r}{S^m} & \frac{n-m+ \frac{1}{2} + |n-m+\frac{1}{2}|}{2(n-m+\frac{1}{2})} \frac{\alpha_n}{\alpha_{n-m}} z^n\end{tikzcd}.  \end{equation}
    Hence we split into two cases to compute the trace.  For $m>n$, \begin{equation}\left(  {S^\dagger}^m S^m - S^m  {S^\dagger}^m \right) z^n = \frac{\alpha_{n+m}}{\alpha_n}z^n.\end{equation}  For $m \leq n$ we have \begin{equation} \left(  {S^\dagger}^m S^m - S^m  {S^\dagger}^m \right) z^n = \left( \frac{\alpha_{n+m}}{\alpha_n}- \frac{\alpha_n}{\alpha_{n-m}}\right) z^n.\end{equation}  Hence we have \begin{equation}\label{series} \tr \left(  {S^\dagger}^m S^m - S^m  {S^\dagger}^m \right)= \displaystyle \sum_{n=0}^{m-1}  \frac{\alpha_{n+m}}{\alpha_n}  + \displaystyle \sum_{n=m}^\infty  \frac{\alpha_{n+m}}{\alpha_n}- \frac{\alpha_n}{\alpha_{n-m}} .\end{equation}  We proceed to show this is absolutely convergent.  Since $\alpha$ is convex, the sequence $\frac{\alpha_{n+1}}{\alpha_n}$ is non-increasing by 4.1 of \cite{zygmund2002trigonometric}, hence the series in \eqref{series} contains all positive terms and the series is absolutely convergent.  In this case the series mostly cancels and we are left with \begin{equation}\displaystyle \lim_{N \rightarrow \infty} \displaystyle \sum_{n=N-m+1}^N \frac{\alpha_{n+m}}{\alpha_n}=m\end{equation} due to our ratio test assumption of \begin{equation}\displaystyle \lim_{n \rightarrow \infty} \frac{\alpha_n}{\alpha_{n+1}}=1.\end{equation} The concave case is similar. Hence we have \begin{equation}\tr \left(  {S^\dagger}^m S^m - S^m {S^\dagger}^m \right)=m.\end{equation} Hence $[{S^\dagger}^m, S^m]$ may be represented by a diagonal matrix with respect to the orthonormal basis $\left\{ \frac{z^n}{\sqrt{\alpha_n}}\right\}_{n=0}^\infty$ of $H_\alpha$.  This is the same orthonormal basis as $\mathcal{D}_0$ hence these Hilbert spaces are equal and the inner products are identical.  This proves this part of the theorem. 
        
    \end{enumerate}  
        
    \end{proof}

\end{theorem}

\section{Generalisation to Weighted Shift}

\begin{definition} Given a Hilbert space $H$ with orthonormal basis $\left\{ z_n\right\}_{n=0}^\infty$ and a weight $r = \left\{ r_n \right\}_{n=0}^\infty$ of complex numbers where $\sup_{n} |r_n| < \infty$, a weighted shift on $H$ is an operator $S_r \in B(H)$ defined by $S_r z_n= r_n z_{n+1}$

\end{definition}

\begin{theorem} Let $\alpha=(\alpha_n)_{n=0}^\infty $ obey  \begin{equation}\displaystyle \lim_{n \rightarrow \infty} \frac{\alpha_{n+1}}{\alpha_{n}}=1\end{equation} and suppose $\alpha$ is concave or convex.  Further suppose $S_r$ is a weighted shift for which \begin{equation}\lim_{n \rightarrow \infty }|r_n| =1.\end{equation}
The commutator of $S_r$ and $S_r^\dagger$ is trace class for all polynomials $f$ and $g$ and we have  \begin{equation}\langle \langle f , g \rangle \rangle_{H_{a_n}}^r = \tr \left( g(S_r)^\dagger f(S_r) - f(S_r) g(S_r)^\dagger \right)=\frac{1}{\pi} \displaystyle \int_{\mathbb{D}} f'(z) \overline{g'(z)} dA(z).\end{equation}  

\begin{proof} We begin by explicitly stating the adjoint shift $S_r^\dagger$ on $f \in H_\alpha$.  We have \begin{equation}S_r^\dagger f = \sum_{n=0}^\infty \frac{\alpha_{n+1}}{\alpha_n} a_{n+1} \overline{r_n} z^n.\end{equation}
We now mimic the proof of the unilateral case; consider $f(z)=g(z)=z^m$.  We have the following: 
\begin{equation}\begin{tikzcd} z^n \arrow[rightarrow]{r}{S_r^m} &  \prod_{i=n}^{n+m-1}r_i z^{n+m}\arrow[rightarrow]{r}{(S_r^\dagger)^m} & \prod_{i=n}^{n+m-1} |r_i|^2 \frac{\alpha_{n+m}}{\alpha_n} z^n\end{tikzcd}.  \end{equation}
Also: \begin{equation}
\begin{aligned}
\begin{tikzcd}[column sep=large]
z^n \arrow[r, "{(S_r^\dagger)^m}"] & 
\displaystyle \prod_{i=n-m}^{n-1} \overline{r_i}
\frac{n-m+ \frac{1}{2} + |n-m+\frac{1}{2}|}{2(n-m+\frac{1}{2})}
\frac{\alpha_n}{\alpha_{n-m}} z^{n-m}
\end{tikzcd}
\\[-1ex]
\begin{tikzcd}[column sep=large]
\arrow[r, "S_r^m"] & 
\displaystyle \prod_{i=n-m}^{n-1}|r_i|^2
\frac{n-m+ \frac{1}{2} + |n-m+\frac{1}{2}|}{2(n-m+\frac{1}{2})}
\frac{\alpha_n}{\alpha_{n-m}} z^n
\end{tikzcd}
\end{aligned}
\end{equation}We split into two cases to compute the trace.  For $m>n$, \begin{equation}\left(  {S_r^\dagger}^m S_r^m - S_r^m  {S_r^\dagger}^m \right) z^n = \prod_{i=n}^{n+m-1}|r_i|^2\frac{\alpha_{n+m}}{\alpha_n}z^n.\end{equation}
For $m \leq n$ we have
\begin{equation} \left(  {S_r^\dagger}^m S_r^m - S_r^m  {S_r^\dagger}^m \right) z^n = \left( \prod_{i=n}^{n+m-1}|r_i|^2\frac{\alpha_{n+m}}{\alpha_n}- \prod_{i=n=m}^{n-1} |r_i|^2\frac{\alpha_n}{\alpha_{n-m}}\right) z^n.\end{equation} Hence we have 
\begin{equation}\begin{split} \tr \left(  {S_r^\dagger}^m S_r^m - S_r^m  {S_r^\dagger}^m \right)= \displaystyle \sum_{n=0}^{m-1} \prod_{i=n}^{n+m-1}|r_i|^2 \frac{\alpha_{n+m}}{\alpha_n}  + \\ \displaystyle \sum_{n=m}^\infty \prod_{i=n}^{n+m-1}|r_i|^2 \frac{\alpha_{n+m}}{\alpha_n}- \prod_{i=n-m}^{n-1} |r_i|^2\frac{\alpha_n}{\alpha_{n-m}}  .\end{split}\end{equation}
By similar arguments we are left with \begin{equation}\lim_{N\rightarrow \infty} \sum_{n=N-m+1}^N \prod_{i=n}^{n+m-1} |r_i|^2 \frac{\alpha_{n+m}}{\alpha_n}.  \end{equation}
By our assumption on the limits of the weights this also uniformly converges to $m$ and the same argument holds on the orthonormal bases.  
\end{proof}
\end{theorem}

\section{Examples}

\begin{definition}[Polylogarithm function] The polylogarithm function defined for $n\in \mathbb{N}$ and $|z|<1$ is given by \begin{equation} \Li_n(z)= \sum_{k=1}^\infty \frac{z^k}{k^n}.  \end{equation}
This is extended to $\mathbb{C}$ by analytic continuation.  We observe that these functions have the property \begin{equation}\Li_{n+1}(z)=\int^z_0 \frac{\Li(t)}{t } dt\end{equation} and $\Li_1(z)=-\ln(1-z)$.  

\end{definition}

\begin{example}
    For the sequences $\alpha_n$ in the top line of the given table, the corresponding $\rkhs$ on the disc has orthonormal basis, reproducing kernel and inner product given on successive lines below $\alpha_n$.  Given $\gamma>-2$ and $\Gamma$ is Euler's gamma function and our $\alpha_n$ may be zero for up to finitely many elements.

\begin{table}[H]
\centering
\scalebox{1}{%}
\begin{tabular}{c c c c c} \toprule
  $\alpha_n$ & $1$  &  $n$ & $n+1$ & $n^2(n-1)$ \\ 
ONB & $\left\{ z^n \right \}_{n=0}^\infty$  &  $\left\{ \frac{1}{\sqrt{n}} z^n \right\}_{n=1}^\infty$ & $\left\{ \frac{1}{\sqrt{n+1}} z^n \right\}_{n=0}^\infty$ & $\left\{ \frac{1}{n\sqrt{n-1}}\right \}_{n=2}^\infty $ \\ 
$K_w (z)$ & $\frac{\overline{w} z }{1-\overline{w} z}$ &  $-\ln (1-\overline{w} z)$ & $-\frac{\ln(1-\overline{w}z)}{\overline{w}z} -1$ & $2 \overline{w}z +  \ln\left((1-\overline{w}z)^{1- \overline{w}z}\right)- \Li_2 (\overline{w}z)$ \\ 
$\langle . , . \rangle$ & $H^2 (\mathbb{D})$  & $\mathcal{D}_0$ & $\mathcal{D}$ &  $\frac{1}{\pi} \displaystyle\int_{\mathbb{D}} f''(z)\overline{g''(z)} dA(z)$ \\ \bottomrule
\end{tabular}%
}
\end{table}

\begin{table}[H]
\centering
\scalebox{1}{%
\begin{tabular}{c c c c} \toprule
  $\alpha_n$ & $\frac{1}{n}$ & $\frac{1}{n+1}$ & $\frac{\Gamma(\beta+1)}{(\gamma+2n+2)^{(1+\beta)}}$ \\  
ONB 
& $\left \{ \sqrt{n} z^n \right \}_{n=1}^\infty$
& $\left\{ \sqrt{n+1} z^n \right\}_{n=0}^\infty$ 
& $\left\{ \frac{(\gamma+2n+2)^{\frac{1+\beta}{2}}}{\sqrt{\Gamma(\beta+1)}}z^n \right\}_{n=0}^\infty$ \\  

$K_w (z)$ 
& $\frac{\overline{w} z}{(1-\overline{w} z)^2}$ 
& $\frac{\overline{w}z(2-\overline{w}z)}{(1-\overline{w}z)^2}$ 
& $\displaystyle \sum_{n=1}^\infty \frac{(\gamma + 2n + 2)^{(1+\beta)}}{(n-k+1)\Gamma (\beta+1)} \overline{w}^n z^n $ \\  

$\langle . , . \rangle$ 
& $\mathcal{A}^2(\mathbb{D})_0$ 
& $\mathcal{A}^2(\mathbb{D})$ 
& $\displaystyle \int_{\mathbb{D}} f(z) \overline{g(z)} |z|^\gamma \left( \log  1 / |z| \right)^\beta \frac{dA(z)}{\pi}$ \\ 
\bottomrule
\end{tabular}%
}
\end{table}

\begin{table}[H]
\centering
\scalebox{1}{%
\begin{tabular}{c c} \toprule
  $\alpha_n$ & $\frac{(n-k+1)\Gamma(\beta+1)}{(\gamma+2n+2)^{(1+\beta)}} \displaystyle\prod_{i=2}^k (n-k+i)^2$ \\  
ONB & $\left\{ \frac{(\gamma + 2n + 2)^{\frac{(1+\beta)}{2}}}{\sqrt{(n-k+1)\Gamma (\beta+1)}} \displaystyle\prod_{i=2}^k \frac{1}{(n-k+i)} z^n \right\}_{n=k}^\infty $ \\  
$K_w (z)$ & $\displaystyle \sum_{n=k}^\infty \frac{(\gamma + 2n + 2)^{(1+\beta)}}{(n-k+1)\Gamma (\beta+1)} \displaystyle\prod_{i=2}^k \frac{1}{(n-k+i)^2} \overline{w}^n z^n $ \\  
$\langle . , . \rangle$ & $\displaystyle \int_{\mathbb{D}} f^{(k)}(z) \overline{g^{(k)}(z)} |z|^\gamma \left( \log  1 / |z| \right)^\beta \frac{dA(z)}{\pi}$ \\ \bottomrule
\end{tabular}%
}
\end{table}

\begin{proof} We show computation of the reproducing kernels and inner products.  Let $\zeta= \overline{w}z$ we have $|\zeta| \leq 1$.
    \begin{itemize}
        \item{$\alpha_n=n^2(n-1)$.  We have \begin{equation}K_w(z)=\displaystyle \sum_{n=2}^\infty \frac{\zeta^n}{n^2(n-1)}.  \end{equation} We differentiate to obtain \begin{equation}\frac{d}{d \zeta} K_w(z)= \displaystyle \sum_{n=2}^\infty \frac{\zeta^{n-1}}{n(n-1)}=\sum_{n=2}^\infty \frac{\zeta^{n-1}}{n-1} -\sum_{n=2}^\infty \frac{\zeta^{n-1}}{n}.  \end{equation} By change of indices we obtain: \begin{equation}\frac{d}{d \zeta} K_w(z)=\sum_{n=1}^\infty \frac{\zeta^n}{n}-\sum_{n=1}^\infty \frac{\zeta^n}{n+1}.  \end{equation} These are computed similarly by differentiating using the geometric series formula, with a change of index for the second term.  We hence obtain \begin{equation} \frac{d}{d \zeta} K_w(z)=-\ln(1-\zeta)+\frac{\ln(1-\zeta)}{\zeta}+1\end{equation} which we integrate by  \cite{gradshteyn2014table} (2.711), (6.254) to obtain \begin{equation}K_w(z)= 2 \zeta + (1- \zeta) \ln(1-\zeta)- \Li_2 (\zeta).
        \end{equation}}
        \item{$\alpha_n= \frac{(n-k+1)\Gamma(\beta+1)}{(\gamma+2n+2)^{(1+\beta)}}\displaystyle\prod_{i=2}^k (n-k+i)^2  $.  We have  \begin{equation} K_w(z)=\frac{1}{\Gamma(\beta+1)}\sum_{n=k}^\infty \frac{\zeta^n(\gamma+2n+2)^{\beta+1}}{n-k+1}\prod_{i=2}^k \frac{1}{(n-k+i)^2} .  \end{equation} This can be computed for any integral values $\beta , k$ for example we consider, $\beta=2$ and $k=3$.  We obtain \begin{equation}\frac{1}{2} \sum_{n=3}^\infty \zeta^n \frac{(\gamma+2n+2)^2}{(n-2)(n-1)^2n^2}.  \end{equation} By partial fractions we obtain: \begin{equation}\begin{split}
        \frac{1}{8} \sum_{n=3}^\infty \zeta^n \bigg( \frac{\gamma^2+12 \gamma +36}{n-2} + \frac{4 \gamma^2 + 16 \gamma}{n-1} -\frac{2\gamma^2 + 8 \gamma + 8}{n^2} - \\ \frac{5 \gamma^2 + 28 \gamma +36}{n} - \frac{4\gamma^2+ 32 \gamma + 64}{(n-1)^2} \bigg) .  
        \end{split}\end{equation} We use standard series formulae results to obtain:
        \begin{equation}\begin{split} \frac{1}{8} \bigg( -(\gamma^2+12\gamma +36)(\zeta^2 \ln(1-\zeta))-(4\gamma^2 + 16 \gamma )(\zeta \ln(1-\zeta) + \zeta^2)+ \\ (2 \gamma^2 + 8 \gamma + 8)\bigg(\zeta+\frac{\zeta^2}{4} - \Li_2 (\zeta)\bigg) + (5 \gamma^2 + 28 \gamma + 36)\bigg(\ln(1-\zeta) + \\ \zeta + \frac{\zeta^2}{2}\bigg) +  (4\gamma^2 + 32 \gamma + 64)(\zeta^2 - \zeta \Li_2 (\zeta))\bigg) .
        \end{split}\end{equation} We see we can calculate these for any values as shown. }
        
    \end{itemize}
\end{proof}
    
\end{example}

\section{Closing Remarks}
The result of our main theorem here gives an instance of the Carey Pincus formula, that is, if $T=X+iY$ is such that $X, Y$ are bounded self-adjoint operators where the commutator $[X, Y]$ is trace class and $T$ acts on a Hilbert space $H$, then for any pair of polynomials \begin{equation} p(x, y)=\sum_{j , k = 1 }^n a_{j k}x^j y^k , q(x, y) = \sum_{j , k = 1}^n b_{j k } x^j y^k , \end{equation} there exists a positive, integrable, compactly supported function $g_T : \mathbb{R}^2 \rightarrow \mathbb{R}$ known as the principal function such that \begin{equation} \tr [p(X , Y) , q(X , Y)]=\frac{1}{2\pi i} \displaystyle \int_{\mathbb{C}} \left( \frac{ \partial p}{\partial x } \frac{\partial q }{\partial y} - \frac{\partial p}{\partial y} \frac{\partial q}{ \partial x}\right) g_T (x , y) dx dy .\end{equation}
Specifically, when $H$ obeys our assumptions, we obtain that $g_T=1$.  Variations of this trace formula are used in context of invariant subspaces of Hilbert space.  Some discussion is found in \cite{zhu2001trace} and \cite{ni2020trace}.  We note that the sequences discussed here are natural to consider; by 4.1 of \cite{katznelson2004introduction} we have that such sequences arise from Fourier transforms of $L^1$ functions.

\end{document}